\documentclass{article}
\usepackage{mathrsfs}
\usepackage{amssymb}
\usepackage{latexsym}
\usepackage{amsmath,amsthm}
\usepackage{amsfonts}
\usepackage{url}

\usepackage{graphicx}

\usepackage{mysymbols}
\usepackage{mytheo}

\newcommand{\F}{F}
\def\eps{\varepsilon}

\renewcommand{\vec}[1]{\ensuremath{\text{{\bf\textrm{#1}}}}}
\newcommand{\prs}{\vec{P}}
\renewcommand{\pr}[1]{\prs\!\tlprn{#1}}
\renewcommand{\mexp}{\vec{E}}
\newcommand{\E}{\mexp}

\renewcommand{\dh}{\operatorname{DH}}
\newcommand{\gkm}{\operatorname{GKM}}

\usepackage{bbm}
\renewcommand{\chr}{\boldsymbol{\mathbbm{1}}} 
\renewcommand{\pred}[1]{\chr_{\left\{ #1 \right\}}}

\newcommand{\setpm}{\set{-1,1}}

\newcommand{\dvc}{d_{\textrm{{\tiny \textup{VC}}}}}
\newcommand{\gor}{\gamma_{\textrm{{\tiny \textup{ORT}}}}}
\newcommand{\poly}{\operatorname{poly}}

\newcommand{\half}{\tfrac{1}{2}}

\begin{document}
\title{VC bounds on the cardinality of nearly orthogonal function classes}
\author{Lee-Ad Gottlieb\thanks{
Weizmann Institute of Science,
\texttt{lee-ad.gottlieb@weizmann.ac.il}
This work was supported in part by The Israel Science Foundation
(grant \#452/08), and by a Minerva grant.
}
, Aryeh Kontorovich\thanks{
Ben Gurion University,
\texttt{karyeh@cs.bgu.ac.il}
}
, Elchanan  Mossel\thanks{
Weizmann institute of Science and U.C.\ Berkeley,
\texttt{mossel@stat.berkeley.edu}. Partially supported by DMS 0548249
(CAREER) award, by ISF grant 1300/08, by a  Minerva Foundation grant
and by an ERC Marie Curie Grant 2008 239317
}
}
\maketitle

\begin{abstract}
We bound the number of nearly orthogonal vectors with fixed VC-dimension over $\setpm^n$.
Our bounds
are of interest
in
machine learning and empirical process theory and improve previous bounds by Haussler.
The bounds are based on a simple projection argument and the generalize to other product spaces.
Along the way we derive tight bounds on the
sum of binomial coefficients in terms of the entropy function.
\end{abstract}

\section{Introduction and statement of results}
The capacity or ``richness'' of a function class $\F$ is a key parameter which makes
a
frequent
appearance in statistics, empirical processes, and machine learning theory
\cite{MR512411,vapnik82,MR757767,talagrand87,MR905334,MR945108,MR1089429,MR1072253}.
It is natural to consider the metric space
$(\F,\rho)$, where $\F\subseteq\setpm^n$ and
\beqn
\label{eq:rhodef}
\rho(x,y) = \oo n\sum_{i=1}^n \pred{x_i\neq y_i}.
\eeqn
A trivial upper bound on the cardinality of $\F$ is
$
2^n$. When $\F$ has
VC-dimension $d$, the celebrated Sauer-Shelah-Vapnik-Chervonenkis lemma \cite{MR0307902} bounds
the cardinality of $\F$ as
\beqn
\label{eq:sauer}
 |\F | \leq \sum_{i=0}^d {n\choose i}.
\eeqn

The notion of cardinality can be refined by considering the packing numbers of the
metric space $(\F,\rho)$.
These are denoted by $M(\eps,d)$, and defined to be the
maximal cardinality of an $\eps$-separated subset of $\F$; in particular
$M(1/n,d) = |\F|$. For general $\eps$, the best packing bound for a maximal
$\eps$-separated subset of $\F$ is due to Haussler \cite{MR1313896}. (A discussion of the
history of this problem may be found therein.) Haussler's upper
bound states that
\beqn
\label{eq:haussler-ub}
M(\eps,d) \leq e(d+1)\paren{\frac{2e}{\eps}}^d.
\eeqn

In this paper, we propose to study the behavior of $M(\eps,d)$ for
$\half-c \le \eps \le \half + c$ (for constant $c$). As explained below,
this corresponds to the case where the vectors of $\F$ are close to
orthogonal. Our interest in this regime stems from applications in machine
learning, where some characterizations and algorithms consider nearly
orthogonal or decorrelated function classes
\cite{DBLP:conf/stoc/BlumFJKMR94,feldman09, dana-david-me-lev/SQ}. Our
main result is Theorem \ref{thm:our-ub} (Section~\ref{sec:upper}), which
sharpens Haussler's estimate of $M(\eps,d)$ as a function of
$d$
and
$\eps\approx\half$.

 \begin{figure}
 \begin{center}
\scalebox{0.6}{
\includegraphics
{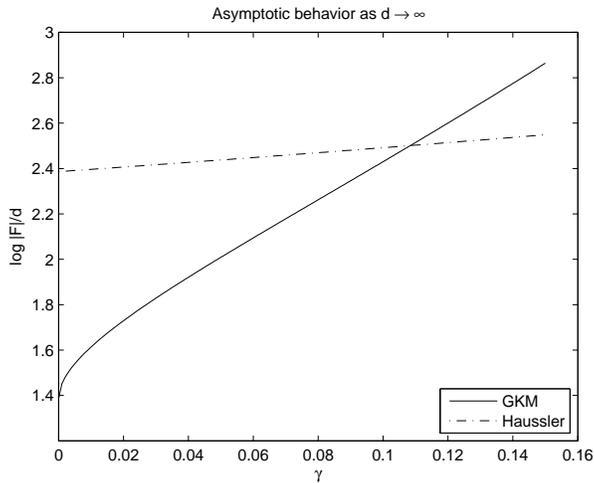}
}
 \end{center}
   \caption{
A comparison of upper bounds.
}
 \label{fig:us-h}
\hide{
g = 0:.001:.15
e=1/2-g/2;
bb = g;
for i=1:length(g)
    bb(i) = gkm_beta(g(i));
end
dh = log(2*exp(1)./e);
us = log(2)*bb;
  x=plot(g,us,'-k'); hold on, y=plot(g,dh,'-.k')
  xlabel('\gamma')
  ylabel('log |F|/d')
  legend([x,y],'GKM','Haussler','Location','SouthEast')
title(['Asymptotic behavior as d \rightarrow \infty'])
}
 \end{figure}

It is convenient to state our results in terms of $\gamma=1-2\eps$ (thus, for
$\eps\approx\half$, we have $\gamma\approx0$).
We will denote
D. Haussler's bound on $|\F|$ in
(\ref{eq:haussler-ub})
by
\beq
M((1-\gamma)/2,d) \leq
\dh(\gamma,d) =
e(d+1)\paren{\frac{4e}{1-\gamma}}^d
\eeq
and
our bound
in Theorem \ref{thm:our-ub}
by
\beq
M((1-\gamma)/2,d) \leq
\gkm(\gamma,d) = 100\cdot 2^{d\beta(\gamma)},
\eeq
where
$\beta:[0,1]\to[2,\infty)$
is defined in (\ref{eq:betadef}).

As $d\to\infty$, our bound asymptotically behaves as
$$ \frac{\ln [\gkm(\gamma,d)]}{d} \to (\ln2) \beta(\gamma) $$
while Haussler's as
$$ \frac{\ln [\dh(\gamma,d)]}{d} \to
\ln \paren{\frac{4e}{1-\gamma}}.
$$
Figure \ref{fig:us-h} gives a visual comparison of these bounds,
illustrating the significant improvement of our bound over Haussler's for
small $\gamma$.

Our analysis has the
additional
advantage of readily extending to $k$-ary
alphabets, while the proof in \cite{MR1313896} appears to be strongly tied
to the binary case.
In Theorem \ref{thm:k-ary} we give what appears to be the first packing bound
for alphabets beyond the binary in terms of (a generalized) VC-dimension (but see
\cite[Lemma 3.3]{alon97scalesensitive}).

We further wish to understand the relationship between $M(\eps,d)$ and
$n$ for fixed $\eps$ and $d$.
It is well known
\cite{RothSeroussi}
that when
$\gamma=1-2\eps=O(1/\sqrt n)$,
we have $M(\eps,d)=O(\poly(n))$.
Since in many cases of interest \cite{MR1965359} the coordinate dimension $n$ may be replaced by its
refinement $\dvc$, it is natural to ask whether a $\poly(n)$ bound on $M(\eps,d)$ is possible for
$\gamma=1-2\eps=O(1/\poly(n))$. We resolve this question in the negative in Theorem \ref{thm:vclb}.


Finally, in Section~\ref{sec:lower} we give a simple improvement of
Haussler's lower bound. Haussler exhibits an infinite family
$\set{\F_n\subseteq\setpm^n}$
for which $\dvc(\F_n)=d$ and
\beqn
\label{eq:haussler-lb}
M(\eps,d) \geq \paren{\oo{2e(\eps+d/n)}}^d.
\eeqn
He notes that the bounds in (\ref{eq:haussler-ub}) and (\ref{eq:haussler-lb}) leave
``a gap from $1/2e$ to $2e$ for the best universal value of the key constant'' and
poses the closure of this gap as an ``intriguing open problem''.
The gap has recently been tightened
to $[1,2e]$
by
Bshouty et al.\ \cite[Theorem 10]{Bshouty2009323}, in a rather general and somewhat involved argument.
Our lower bound
in Theorem \ref{thm:c0}
achieves the same tightening via a much simpler construction.

\section{Definitions and notation}
\label{sec:def}
Our basic object is the metric space
$(\F,\rho)$, with
$\F\subseteq\setpm^n$ and
the normalized Hamming distance $\rho$ defined in (\ref{eq:rhodef}).
The inner product
\beq
\iprod{x}{y}:=\oo n\sum_{i=1}^n x_i y_i,
\qquad
x,y\in\F
\eeq
endows $\F$ with Euclidean structure. The distance and inner product have a simple relationship:
\beqn
\label{eq:gamrho}
2\rho(x,y) + \iprod{x}{y} = 1.
\eeqn
We denote the natural numbers by $\N=\set{1,2,\ldots}$, and
for $n\in\N$, we write $[n]=\set{0,1,\ldots,n-1}$.
For $I=(i_1,i_2,\ldots,i_k)\subseteq
[n]
$, we denote the projection of $\F$ onto $I$ by
\beqn
\label{eq:projdef}
\evalat{\F}{I}
= \set{(x_{i_1},\ldots,x_{i_k}): x\in\F}\subseteq\setpm^k.
\eeqn
We say that $\F$ {\em shatters} $I$ if
$\evalat{\F}{I} =\setpm^k$
and define the Vapnik-Chervonenkis dimension of $\F$ to be the cardinality
of the largest shattered index sequence $I$:
\beq
\dvc(\F) = \max\set{|I| : I\subset
[n]
,
\evalat{\F}{I} =\setpm^k}.
\eeq

We define
$\gamma=\gor(\F)$
by
\beqn
\label{eq:gorth}
 \gor(\F) = \max\set{|\iprod{x}{y}|:x\neq y\in\F
}.
\eeqn
In words, $\gor(\F)$ is the smallest $\gamma\geq0$ such that
all distinct pairs $x,y\in\F$ are
``orthogonal to accuracy $\gamma$''.
Whenever (\ref{eq:gorth}) holds for some $\gamma$, we say that $\F$ is $\gamma$-orthogonal.
\hide{
For $0\leq\gamma\leq 1$
we denote by
$J(\gamma,\F)$  the cardinality of the largest $\gamma$-orthogonal
subset of $\F$:
\beq
J(\gamma,\F) := \max\set{ |\F'| : \F'\subset\F,  \gor(F')\leq\gamma}.
\eeq
}

We will use $\ln$ to denote the natural logarithm and $\log\equiv\log_2$.

\section{Upper estimates on nearly orthogonal sets}\label{sec:upper}

\subsection{Preliminaries: entropy and $\beta$}
Recall the binary entropy function, defined as
\beqn
\label{eq:hdef}
H(x)=-x \log x - (1-x) \log (1-x).
\eeqn
In the range $[0,1]$, this function is symmetric about $x=\half$,
where it achieves its maximum value of 1.

\begin{figure}
\hide{
g = 0:.0001:.25
bb = g;
for i=1:length(g)
    bb(i) = gkm_beta(g(i));
end
plot(g,bb,'.')
xlabel('\gamma')
ylabel('\beta(\gamma)')
}
 \begin{center}
\scalebox{0.6}{
\includegraphics
{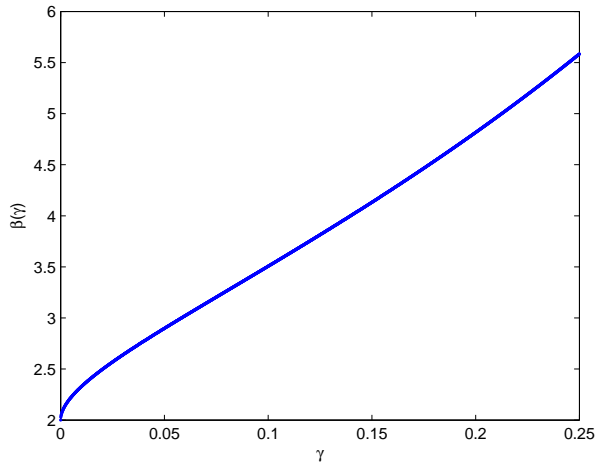}
}
 \end{center}
   \caption{
The function $\beta(\gamma)$.
}
 \label{fig:beta}
 \end{figure}

Since $H$ is increasing on $[0,\half]$, it has a well-defined inverse
on this domain, which we will denote by $H\inv:[0,1]\to[0,\half]$.
We define the function $\beta:[0,1]\to[2,\infty)$ by
\beqn
\label{eq:betadef}
\beta(\gamma) = \oo{ H\inv[
\log(2/(1+\gamma))
]}.
\eeqn
Figure \ref{fig:beta}
illustrates the behavior of $\beta$ on $[0,\tfrac{1}{4}]$.

A sharp bound on $\sum_{i=0}^d {n\choose i}$ in terms of $H$ is given
in Lemma \ref{lem:binom}.

\subsection{Main result}

\bethn
\label{thm:our-ub}
Let $\F\subseteq\setpm^n$ with $
1\leq
d=\dvc(\F)
\leq n/2
$ and $\gamma=\gor(\F)$.
Then
\beq
\abs{\F} &\leq &
100\cdot
2^{d\beta(\gamma)}
\eeq
where $\beta(\cdot)$ is defined in (\ref{eq:betadef}).
\enthn
\bepf
Let $r<n$ be unspecified for the moment
and
choose $I\subset[n]$, $|I|=r$ uniformly at random.
Define $\pi=\pi_I$ to be the coordinate projection
of $F$ onto $I$ as defined in (\ref{eq:projdef}).
Let $x$ and $y$ be two uniformly random elements of $\F$, and let $A$ be the event that
$\pi(x)=\pi(y)$;
thus, $P(A)$ is the probability that $x$ and
$y$ are mapped to the same vector.
The latter
is upper-bounded by the sum of the probability that
$x$ and $y$ are the same vector,
and the
probability that $x$ and $y$ are distinct vectors but are mapped to the same vector.
The first event occurs with probability exactly
$
|\F|\inv
$. We claim that
the second event occurs with probability less
than $(\frac{1}{2} + \frac{1}{2}\gamma)^r$.
To see this, suppose that
the two vectors
$x,y$ agree on $\eta$ fraction of the coordinates.
Then $\eta \leq \half+\half \gamma$ and the probability that they agree on
one random coordinate is exactly $\eta$.
The probability they agree on two
coordinates
is $\eta (n \eta - 1)/(n-1)$, and so forth.
Thus, the probability that they agree on
$r$ coordinates is
$$
\eta  (n \eta - 1)/(n-1)\cdot \ldots \cdot (n \eta -(r-1))/(n-(r-1)) <
\eta^r \leq (\half+\half\gamma)^r.$$
By the union bound, we have
\beqn
\label{eq:PAub}
P(A) <
|\F|\inv
+
\paren{\oo2+\oo2\gamma}^r.
\eeqn
As a lower bound on $P(A)$, we claim
\beqn
\label{eq:PAinv}
P(A)\inv
\leq
\sum_{i=0}^d \binom{r}{i}
.
\eeqn
Indeed,
if
$E$ is any
finite set
equipped with distribution
$P_E$, then the probability of collision (i.e., drawing $e,e'\in E$ independently according to $P_E$
and having $e=e'$) is given by
$P_E(e=e')=\sum_{e\in E} P_E(e)^2$.
Now by Jensen's inequality,
$$ |E|^{-2} = \paren{\sum_{e\in E} |E|\inv P_E(e)}^2
\leq \sum_{e\in E} |E|\inv P_E(e)^2,$$
which implies
\beqn
\label{eq:P_E}
P_E(e=e')=\sum_{e\in E} P_E(e)^2\geq |E|\inv.
\eeqn

Let us
denote
the event that $\pi(x)=\pi(y)$ conditioned on $I$
by $A\gn I$,
and
write
$P_\pi$
for
the distribution on $\F':=\evalat{\F}{I}$ induced by $\pi$.
Then
we have
\beq
P(A\gn I) &=& \sum_{x'\in\F'} P_\pi(x')^2 \\
&\ge&  |\F'|\inv \\
&\ge& \paren{\sum_{i=0}^d \binom{r}{i}}\inv
,
\eeq
where
the first inequality is seen by taking $E=F'$ and $P_E=P_\pi$ in (\ref{eq:P_E})
and the second holds by Sauer's Lemma (\ref{eq:sauer}).
The claim (\ref{eq:PAinv})
follows by
averaging over all the $I$s.

\hide{
This lower bound is achieved by the uniform distribution $P_E(x)\equiv|E|\inv$, and therefore
\beqn
\label{eq:PAlb}
P(A) \geq |\F_r|\inv
\eeqn
where $\F_r$ is the restriction of $\F$ to the $r$ coordinates we projected on.
By Sauer's Lemma and Lemma \ref{lem:binom}, we have
\beqn
\label{eq:Fgrbd} 
|\F_r| \leq
\sum_{i=0}^d \binom{r}{i} < .9964 \cdot 2^{rH(d/r)},
\eeqn
which holds for all $2\leq d\leq r/2$.
}

Combining
(\ref{eq:PAub})
and
(\ref{eq:PAinv})
with Lemma \ref{lem:binom},
we get the key inequality
\beqn
\label{eq:key}
1.02 \cdot
2^{-rH(d/r)}
<
\frac{1}{|\F|} +
\paren{\oo2+\oo2\gamma}^r,
\eeqn
valid
for all
integer
$r\in[2d,n]$.
We
choose the value
\beq
r^* = \ceil{\beta(\gamma)d}
\eeq
where the function $\beta(\cdot)$ is defined in (\ref{eq:betadef}).
It is straightforward to verify from the definition of $\beta(\cdot)$ that
for this choice of $r^*$, we have
\beq
2^{-r^*H(d/r^*)} \geq \paren{\oo2+\oo2\gamma}^{r^*}
\eeq
and therefore
\beq
\label{eq:Frbd}
.02 \cdot
2^{-r^*
}
<
|\F|\inv
;
\eeq
combining this with (\ref{eq:key}) yields
\beq
|\F| &\leq & 50 \cdot 2^{\ceil{\beta(\gamma)d}
} \\
&\leq& 100 \cdot 2^{\beta(\gamma)d } .
\hide{
\cdot
2^{
-\paren{1+
d\inv
}\log(\half+\tfrac{\gamma}{2})
}\\
&=& 277
\paren{\frac{2}{1+\gamma}}^{1+d\inv}
2^{\beta(\gamma)d} \\
&\leq&
784 \cdot 2^{\beta(\gamma)d}.
}
\eeq
\enpf

\section{Generalization to $k$-ary alphabets}
\label{sec:k-ary}
Here we extend our upper bound analysis to $k$-ary ($k\ge3$) alphabets.
First, we must
generalize the notion of orthogonality.
Since two vectors $x,y$ drawn uniformly from
$[k]^n
$
agree in expectation on $n/k$ coordinates, we may define
$\gamma_k(x,y)$ by
\beqn
\label{eq:gamrhok}
\frac{k}{k-1}\rho(x,y) + \gamma_k(x,y) = 1,
\eeqn
where $\rho$ is
the normalized Hamming distance
defined in (\ref{eq:rhodef}).
Analogously, we
define
$\gor^k(\F)$ by
\beqn
\label{eq:gork}
 \gor^k(\F) = \max\set{|\gamma_k(x,y)|:x\neq y\in\F
}.
\eeqn
\hide{
$J_k(\gamma,\F)$ by
\beq
J_k(\gamma,\F) := \max\set{ |\F'| : \F'\subset\F,
\gamma_k(x,y)\leq \gamma
~\mbox{for all}~
x\neq y\in\F'
}.
\eeq
}
The notion of VC-dimension
has various generalizations
to $k$-ary alphabets
\cite{DBLP:conf/focs/Haussler89,DBLP:journals/ml/Natarajan89,pollard84,MR1089429}.
Among these, we consider Pollard's P(seudo)-dimension, Natarajan's G(raph)-dimension,
and the GP-dimension; these are defined in equations (13,14,15) of \cite{209962},
respectively.
In the sequel we continue to write $\dvc(\F)$ to denote one of these
combinatorial dimensions, without specifying which one we have in mind.
This convention is justified
by a
common generalized
Sauer's Lemma shared by these three quantities, due to
Haussler and Long \cite[Corollary 3]{209962}:
\beqn
\label{eq:HL}
|\F| \leq \sum_{i=0}^{\dvc(\F)} \binom{n}{i}k^i
.
\eeqn
A sharp bound on the rhs of (\ref{eq:HL})
is given in Lemma \ref{lem:binom-k}.

Our main result is readily generalized to $k$-ary alphabets:

\bethn
\label{thm:k-ary}
Let $\F\subseteq[k]^n$
with $
\frac{6k}{k+1.6}\leq
d=\dvc(\F)
\leq \frac{nk}{k+1.6}
$ and $\gamma=\gor^k(\F)$.
Then
\beq
\abs{\F} &\leq &
34 k^d 2^{d/\delta(\gamma,k)}
\eeq
where
$\delta(\gamma,k)$
is the largest $x\in[0,k/(k+1)]$ for which
$ x\log k + H(x) \leq \log (k/(1+(k-1)\gamma))$
holds.
\enthn
\noindent {\bf Remark:}
The function $\delta:(0,1)\times\N\to(0,1)$ is readily computed numerically.
\bepf
Repeating the argument in Theorem \ref{thm:our-ub} (with
the generalized Sauer Lemma (\ref{eq:HL})),
we have
\beq
\paren{\sum_{i=0}^{d} \binom{r}{i}k^i}\inv
<
|\F|\inv + \paren{\oo k + \frac{k-1}{k}\gamma}^r.
\eeq
Applying the bound in Lemma \ref{lem:binom-k}, we have that for
$
\frac{6k}{k+1.6}\leq
d
\leq \frac{rk}{k+1.6}
$,
\beq
1.06 \cdot 2^{-rH(d/r) - d \log k}
<
|\F|\inv + \paren{\oo k + \frac{k-1}{k}\gamma}^r.
\eeq

Now we seek the minimum integer
$r\in[\frac{k+1.6}{k}d,n]$ that ensures
$$ d\log k + rH(d/r) \le r\log
(k/(1+(k-1)\gamma))
. $$
To this end, we consider the following
inequality
in $x$
\beqn
\label{eq:AHB}
 x\log k + H(x) \leq \log (k/(1+(k-1)\gamma))
.
\eeqn
Note that the inequality (\ref{eq:AHB}) is satisfied at $x=0$ and define
$x^*\equiv\delta(\gamma,k)$ to be the
largest $x\in[0,k/(k+1.6)]$ satisfying it (the proof of
Lemma \ref{lem:binom-k} shows that
the lhs of (\ref{eq:AHB}) is monotonically increasing in this range).
Taking $r^* = \ceil{d/x^*}$, we have
\beq
.06 \cdot 2^{-r^*H(d/r^*) - d \log k} < |F|\inv,
\eeq
which rearranges to
\beq
|F| &<& 17\cdot 2^{r^*H(d/r^*) + d \log k} \\
&\leq& 34 k^d  2^{d/\delta(\gamma,k)},
\eeq
as claimed.
\hide{

Note that for $x=0$, we have lhs(\ref{eq:AHB}) $<$ rhs(\ref{eq:AHB}),
and for $x=1$, lhs(\ref{eq:AHB}) $>$ rhs(\ref{eq:AHB}). Since all the
terms in (\ref{eq:AHB}) are continuous in $x$, we may define $\delta(\gamma,k)$
as the smallest $x>0$ for which
$$ x\log k + H(x) = \log (k/(1+(k-1)\gamma))$$
holds.

Taking $r=\ceil{d/\delta(\gamma,k)}$ yields the claimed bound.
}
\enpf
\section{Polynomial upper bounds for small $\gamma$}\label{sec:polynomial}
The bounds of Haussler (\ref{eq:haussler-ub})
and Theorem \ref{thm:our-ub} obscure the dependence of
$|\F|$ on its coordinate dimension $n$.
It is well known that when
$\gor(\F)=O(1/\sqrt n)$,
we have $
\F
=O(\poly(n))$.
(In the degenerate case $\gor(\F)=0$, linear algebra
gives $|\F|\leq n+1$.)

Roth and Seroussi
\cite{RothSeroussi}
developed a powerful technique for
bounding $|\F|$ in terms of $n$ and $\gamma$.
Let $0<\rho_{\min} \leq \rho_{\max}$ be such that
\beq
\rho_{\min} \leq n\rho(x,y) \leq \rho_{\max}
\eeq
for all $x,y\in\F$.
Then
\cite[Proposition 4.1]{RothSeroussi}
shows that
\beq
1-|\F|\inv
\le \paren{1-\frac{1}{n}}
\paren{\frac{
\rho_a}{\rho_g}}^2
\eeq
where $\rho_a=\tfrac{1}{2}(\rho_{\min}+\rho_{\max})$
and
$\rho_g=\sqrt{\rho_{\min}\rho_{\max}}$.
Recalling the relation in (\ref{eq:gamrho}), we have
\beq
\rho_{\max}=\frac{n}{2}(1+\gamma),~
\rho_{\min}=\frac{n}{2}(1-\gamma),~
\rho_a=\frac{n}{2},~
\rho_g=\frac{n}{2}\sqrt{1-\gamma^2},
\eeq
which implies the following bound on $|\F|$:
\beq
1-|\F|\inv
\leq
\paren{1-\frac{1}{n}}\oo{1-\gamma^2}.
\eeq
Note that
when $\gamma^2\geq n^{-1}$, the right-hand side is least 1 and the bound is rendered vacuous; thus
the nontrivial regime is $\gamma^2 < n\inv$. In particular, taking $\gamma=1/(c\sqrt n)$ for $c>1$ yields the bound
\beqn
\label{eq:cn}
|\F| \leq \frac{c^2n-1}{c^2-1}.
\eeqn

Since in many situations, the VC-dimension $\dvc$ is a refinement of the coordinate dimension $n$,
it is natural to ask if a bound similar to (\ref{eq:cn}) holds with $\dvc$ in place of $n$.
We resolve this question strongly in the negative:

\bethn
\label{thm:vclb}
Let $a>0$ be some constant.
Then there
infinitely many $n\in\N$ for which
there
is an $\F\subseteq\setpm^n$
such that
\bit
\item[(a)] $\gamma= d^{-a}
$
\item[(b)] $\abs{\F}
=
\floor{\exp\paren{
cn^ {\frac{1}{2a+1}}
}}
$
\eit
where $\gamma=\gor(\F)$ , $d=\dvc(\F)$ and $c
$ is an absolute constant.
\enthn
\bepf
Let $\F$ be an $m\times n$ matrix whose entries are independent
symmetric Bernoulli $\setpm$ random variables; we shall identify the rows
of $\F$ with the functions in $\F$. Then for $f,g\in\F$, we have
$$ \E\iprod{f}{g}=0$$
and by Chernoff's bound
$$ \pr{ \abs{\iprod{f}{g}} > \gamma} \leq 2\exp(-n\gamma^2/2)$$
for all
$n\in\N$ and
$\gamma>0$.
The union bound implies that for $n$ large enough
there exists an
$\F\subseteq\setpm^n$
with
$\gor(\F)\leq\gamma$
and
$$
|\F|=\floor{\exp(n\gamma^2/4)}
.$$
The claim follows from
the relation
$$d= \dvc(\F) \leq \log_2|\F|  \leq
n\gamma^2/4\ln2
$$
and our choice of
$$ \gamma=d^{-a}.$$
\hide{
we establish the claim.
$$ d^{2a+1} \leq cn $$
and
$$ \gamma^2 \leq d^{-2a} \leq cn^{-\frac{2a}{2a+1}}.$$
Therefore,
\beq
m &=& C\exp(n\gamma^2) \\
&\geq& C\exp(c n^ {-\frac{1}{2a+1}}) \\
&\geq& C \exp(c d)
\eeq
as claimed.
}

\enpf

An alternative estimate may be obtained via the Gilbert-Varshamov bound \cite{Gilbert52,Varshamov57}.

\section{A lower bound on the universal constant $c_0$}\label{sec:lower}
Haussler's upper
(\ref{eq:haussler-ub}) and lower (\ref{eq:haussler-lb}) bounds imply
the existence of a
universal $c_0$ for which
the packing number
$M(\eps,d)$ grows as
$\Theta((c_0/\eps)^d)$
in $\eps$
for constant $d$.
More precisely,
\bit
\item[(i)] $M(\eps,d)=O(d(c_0/\eps)^d)$
for all $n,\F\subseteq\setpm^n$ with $\dvc(\F)= d$
\item[(ii)] $M(\eps_n,d)=\Omega((c_0/\eps_n)^{d_n})$
for some infinite family $(\eps_n,d_n,\F_n\subseteq\setpm^n)$ with $\dvc(\F_n)= {d_n}$.
\eit
The bounds in (\ref{eq:haussler-ub}, \ref{eq:haussler-lb}) peg $c_0$ at
$1/2e \leq c_0 \leq 2e$.
An improved lower bound of $c_0\geq1$
may be obtained
essentially ``for free'' (cf.
\cite[Theorem 10]{Bshouty2009323}):
\bethn
\label{thm:c0}
There exists an infinite family
$(\eps_n,d_n,\F_n\subseteq\setpm^n)$
for which
\bit
\item[(a)] $\dvc(\F_n)= {d_n}$
\item[(b)] $M(\eps_n,d)=(1/\eps_n)^{d_n}$
\eit
\enthn
\bepf
For $n=1,2,\ldots$, put $\eps_n=\half$, $d_n=n$, and $\F_n\subset\setpm^n$ to be the rows of $H_{2^n}$,
the Hadamard matrix of order $2^n$. The latter may be defined recursively via
\beq
H_1 = [1]
\eeq
and
\beq
H_{2^{n+1}} =
\sqprn{
\begin{array}{rr}
H_{2^n} & H_{2^n} \\
H_{2^n} & -H_{2^n}
\end{array}
}.
\eeq
It is well known (and elementary to verify) that $\dvc(\F_n)=n$ and
that $\gor(\F_n)=0$. Thus $\F_n$
is a $\half$-separated set of size $2^n$.
\enpf

\section{Technical Lemmata}
\label{sec:tech}

Our main result
in Theorem \ref{thm:our-ub}
requires a sharp
estimate on the sum of the binomial coefficients.
It is well known
\cite{1121738}
that for $d \le \frac{n}{2}$,
$\sum_{i=0}^d \binom{n}{i} \le 2^{nH(d/n)}$,
but we need to obtain a slightly tighter bound.
\belen
\label{lem:binom}
For $1 \le d \le \frac{n}{2}$, we have
\beq
\sum_{i=0}^d \binom{n}{i} &<& \delta \cdot 2^{nH(d/n)},
\eeq
where $\delta = 0.98$.
\enlen

\noindent {\bf Remark:}
The bound $\delta$ can be further tightened, at the expense of a more complicated proof.
Note however that when $d=n/2$ the summation is equal to $\frac{1}{2} 2^{nH(d/n)}$, so $\delta$ cannot be
taken as a constant better than $\half$.

\bepf
Recall Stirling's approximation
$i! = \sqrt{2\pi i}\paren{\frac{i}{e}}^i e^{\lambda_i}$ where
$\frac{1}{12i+1} < \lambda_i < \frac{1}{12i}$.
Also note that for $0 \le i \le n$,
$$\frac{1}{12n}-\frac{1}{12(n-i)+1}-\frac{1}{12i+1}
=\frac{-144n^2+122ni-144i^2-12n}{(12n)(12n-12i+1)(12i+1)}
\le 0.$$
Thus,
\beq
\binom{n}{i}
& = & \frac{n!}{i!(n-i)!}   \\
& \le &
e^{\frac{1}{12n}-\frac{1}{12(n-i)+1}-\frac{1}{12i+1}} \cdot
\sqrt{\frac{n}{2\pi i(n-i)}} \cdot
\frac{n^n}{i^i (n-i)^{n-i}}     \\
& <  &
\frac{1}{\sqrt{2\pi i (1-i/n)}} \cdot
(i/n)^{-i} (1-i/n)^{-(n-i)} \\
& =  &
\frac{1}{\sqrt{2\pi i (1-i/n)}} \cdot 2^{nH(i/n)}.
\eeq

We first prove Lemma~\ref{lem:binom} for small values of $d$, in particular
$1 \le d < {n}/{4}$.
Note that for $i\le d < n/4$ we have
\beq
\binom{n}{i-1} = \frac{i}{n-i+1} \binom{n}{i} < \frac{1}{3} \binom{n}{i},
\eeq
and therefore
\beq
\sum_{i=0}^d \binom{n}{i}
& < & 1.5 \binom{n}{d}	\\
& < & \frac{1.5}{\sqrt{2\pi d (1-d/n)}} \cdot 2^{nH(d/n)} 	\\
& < & .7 \cdot 2^{nH(d/n)}.
\eeq

We now turn to the case of large $d$, that is $\frac{n}{4} \le d \le \frac{n}{2}$.
If $\sum_{i=0}^d \binom{n}{i} < 0.5 \cdot 2^{nH(d/n)}$, then Lemma~\ref{lem:binom}
immediately holds,
so we may assume that $Z := \sum_{i=0}^d \binom{n}{i} \geq 0.5 \cdot 2^{nH(d/n)}$.
We will show that in this case, much of the weight of the sum is distributed among at least 
$\Omega(\sqrt{n})$ coefficients. We will use this fact in conjunction with the standard entropy argument, see e.g.,  \cite{1121738} to obtain the desired result.

Now, we have for all $i \leq d$ (when $\frac{n}{4} \le d \le \frac{n}{2}$),
\beq
\binom{n}{i}
\leq  \binom{n}{d}
<  \frac{1}{\sqrt{2\pi d (1-d/n)}} \cdot 2^{nH(d/n)}
<  \frac{2}{\sqrt{\pi n}} 2^{nH(d/n)}
\leq  \frac{4 Z}{\sqrt{\pi n}}
\eeq

Consider the random vector $(X_1,\ldots,X_n)$ uniformly distributed in
$\{ x : \{0,1\}^n : \sum_i x_i \leq d \}$. Then for all $0 \leq r \leq d$ we have: 
\[
P\sqprn{\sum_{i=1}^n X_i = r} = Z^{-1} \binom{n}{i} \leq \frac{4}{\sqrt{\pi n}},
\]
and therefore 
\[
P\sqprn{\sum_{i=1}^n X_i \geq d - \frac{\sqrt{\pi n}}{8} + 1} \leq  \frac{\sqrt{\pi n}}{8} \frac{4}{\sqrt{\pi n}} \leq 
\frac{1}{2},
\]
which implies
\beq
\E\sqprn{\sum_{i=1}^n X_i} &\leq& d P\sqprn{\sum_{i=1}^n X_i \geq d - \frac{\sqrt{\pi n}}{8} + 1} + 
\paren{d - \frac{\sqrt{\pi n}}{8}} \paren{1 - P\sqprn{\sum_{i=1}^n X_i \geq d - \frac{\sqrt{\pi n}}{8} + 1}} \\
&\leq& \frac{1}{2}\paren{d + d - \frac{\sqrt{\pi n}}{8}} = d - \frac{\sqrt{\pi n}}{16}.
\eeq

Hence, we obtain
\beq
H(X_1,\ldots,X_n)
& \le & n H(X_i) = n H(\E[X_i]) \\
& < & n H \left( \frac{d}{n}
-\frac{\sqrt{\pi}}{16\sqrt{n}} \right) 	\\
& = & n H \left( \frac{d}{n} \right)
- n \left( H \left( \frac{d}{n} \right)
- H \left( \frac{d}{n} -\frac{\sqrt{\pi}}{16\sqrt{n}} \right) \right) 	\\
& < & n H \left (\frac{d}{n} \right)
- n \left( H \left( \frac{1}{2} \right) - H \left( \frac{1}{2}-\frac{\sqrt{\pi}}{16\sqrt{n}} \right) \right),
\eeq
where the second inequality uses the monotonicity of the binary entropy function $H$ at
$[0,\half]$,
and the third uses the concavity of $H$. Noting that the Taylor series expansion of
$H(x)$ around $\frac{1}{2}$ is equal to
$1 - \frac{1}{2 \ln 2} \sum_{j=1}^\infty \frac{(1-2x)^{2j}}{j(2j-1)}
< 1 - \frac{(1-2x)^{2}}{2 \ln 2}$, we have that
\beq
H \left( \frac{1}{2} \right)
- H \left( \frac{1}{2}-\frac{\sqrt{\pi}}{16\sqrt{n}} \right)
> \frac{1}{2 \ln 2} \frac{\pi}{64 n},
\eeq
from which we conclude that
\beq
H(X_1,\ldots,X_n) < n H(d/n) - \frac{\pi}{128 \ln 2}.
\eeq
Hence, we have
\beq
\sum_{i=0}^d \binom{n}{i}
& = & 2^{H(X_1,\ldots,X_n)} \\
& < & 2^{-\frac{\pi}{128 \ln 2}} 2^{n H(d/n)}	\\
& < & .98 \cdot 2^{n H(d/n)},
\eeq
where the first identity holds because $H(Y)=\log|\mathrm{supp}(Y)|$ when $Y$ is uniformly distributed on its support.
This completes the proof.
\enpf

\hide{
\section{Technical Lemmata}
\label{sec:tech}
We will need a simple monotonicity property of
the entropy function
$H$
defined in (\ref{eq:hdef}):
\belen
The function $F:\N\to\R$ given by
$F(n)=nH(x/n)$ is nondecreasing in $n$ for all $x\in[0,1]$.
\enlen
\bepf


The derivative
$$ F'(n) = \log(n/(n-x)) $$
is nonnegative for all $x\in[0,1]$.
\enpf

Our main result
in Theorem \ref{thm:our-ub}
requires a sharp
estimate on the sum of the binomial coefficients.
It is well known
\cite{1121738}
that for $d \le \frac{n}{2}$,
$\sum_{i=0}^d \binom{n}{i} \le 2^{nH(d/n)}$,
but we need a slightly tighter bound.

\belen
\label{lem:binom}
For $2 \le d \le \frac{n}{2}$, we have
\beq
\sum_{i=0}^d \binom{n}{i} &<& .9964 \cdot 2^{nH(d/n)}.
\eeq
\enlen

\bepf
Recall Stirling's approximation
$i! = \sqrt{2\pi i}\paren{\frac{i}{e}}^i e^{\lambda_i}$ where
$\frac{1}{12i+1} < \lambda_i < \frac{1}{12i}$.
Also note that for $0 \le i \le n$,
$$\frac{1}{12n}-\frac{1}{12(n-i)+1}-\frac{1}{12i+1}
=\frac{-144n^2+122ni-144i^2-12n}{(12n)(12n-12i+1)(12i+1)}
\le 0.$$
Thus,
\beq
\binom{n}{i}
& = & \frac{n!}{i!(n-i)!}   \\
& \le &
e^{\frac{1}{12n}-\frac{1}{12(n-i)+1}-\frac{1}{12i+1}} \cdot
\sqrt{\frac{n}{2\pi i(n-i)}} \cdot
\frac{n^n}{i^i (n-i)^{n-i}}     \\
& <  &
\frac{1}{\sqrt{2\pi i (1-i/n)}} \cdot
(i/n)^{-i} (1-i/n)^{-(n-i)} \\
& =  &
\frac{1}{\sqrt{2\pi i (1-i/n)}} \cdot 2^{nH(i/n)}.
\eeq

This gives us an upper bound on the summation of binomial terms
$$\sum_{i=0}^d {\binom{n}{i}}
= 1+\sum_{i=1}^d {\binom{n}{i}}
< 1+\sum_{i=1}^d
\frac{1}{\sqrt{2\pi i(1-i/n)}} \cdot
2^{nH(i/n)}.$$
It is left to upper bound the summation of the right hand side. We first note
that the function $\frac{2^{nH(i/n)}}{\sqrt{2\pi i(1-i/n)}}$ is monotone
increasing in the range $1 \le i < \frac{n}{2}$, and so the terms of this
summation are strictly increasing. To prove this, define the continuous function
\beq
f(x)
& = & \ln \left( \frac{1}{\sqrt{2\pi x (1-x/n)}} \cdot 2^{nH(x/n)} \right) \\
& = & \ln \frac{1}{\sqrt{2\pi}} - \frac{1}{2} \ln x - \frac{1}{2} \ln (1-x/n)
  -x\ln (x/n) - (n-x)\ln(1-x/n)
\eeq
We show that $f(x)$ is strictly increasing in the relevant range
$1 \le x < \frac{n}{2}$. The first and second derivatives of $f(x)$ are
\beq
f'(x)  & = & - \frac{1}{2x} + \frac{1}{2(n-x)} - \ln x + \ln (n-x), \\
f''(x) & = & \frac{1}{2x^2} + \frac{1}{2(n-x)^2} - \frac{1}{x} - \frac{1}{n-x}.
\eeq
Now, $f''(x)$ is negative in the relevant range, since
$\frac{1}{2x^2}- \frac{1}{x} = \frac{1-2x}{2x^2} < 0$ and
$\frac{1}{2(n-x)^2} - \frac{1}{n-x} = \frac{1-2n + 2x}{2(n-x)^2} \le 0$. It follows that
$f'(x)$ is decreasing in this range, and since $f'(n/2) = 0$, $f'(x)$ is always positive
in the range. It follows that $f(x)$ is strictly increasing in the relevant range, and
so the terms of the summation are increasing.

Henceforth, we assume that $n \ge 41$. (The finite number of cases not covered by our
analysis are confirmed numerically.) Let $p=\frac{n}{2} - \sqrt{n}/c$, where
$c=2.8 < 2 \sqrt{2}$. First note that the derivative of $2^{nH(i/n)}$ with respect to $i$
is equal to $2^{nH(i/n)} \ln (\frac{n}{i} - 1) di$. Below, we will demonstrate that for
$1 \le i \le p$ (and $n \ge 41$), we have that
\beq
\ln (\frac{n}{i} - 1) > \frac{1}{\sqrt{i}}.
\eeq
It follows that for $d \le p$,
\beq
\sum_{i=1}^d {\binom{n}{i}}
& < & \sum_{i=1}^d \frac{1}{\sqrt{2\pi i(1-i/n)}} \cdot 2^{nH(i/n)} \\
& = & \frac{1}{\sqrt{\pi}} \left[
      \frac{1}{\sqrt{2 d(1-d/n)}} \cdot 2^{dH(d/n)} +
      \sum_{i=1}^{d-1} \frac{1}{\sqrt{2 i(1-i/n)}} \cdot 2^{nH(i/n)} \right] \\
& = & \frac{1}{\sqrt{\pi}} \left[
      \frac{1}{\sqrt{2 d(1-d/n)}} \cdot 2^{dH(d/n)} +
      \int_{1}^d \frac{1}{\sqrt{2 i(1-i/n)}} \cdot 2^{nH(i/n)} di \right] \\
& < & \frac{1}{\sqrt{\pi}} \left[
      \frac{1}{\sqrt{d}} \cdot 2^{dH(d/n)} +
      \int_{1}^d \frac{1}{\sqrt{i}} \cdot 2^{nH(i/n)} di \right] \\
& < & \frac{1}{\sqrt{\pi}} \left[
      \frac{1}{\sqrt{d}} \cdot 2^{nH(d/n)} +
      \int_1^d \ln (\frac{n}{i} - 1) \cdot 2^{nH(i/n)} di \right] \\
& < & \frac{1}{\sqrt{\pi}} \left( \frac{1}{\sqrt{d}}+1 \right) \cdot 2^{nH(d/n)} \\
& < & .97 \cdot 2^{nH(d/n)},
\eeq
where the final inequality follows from $d \ge 2$. Hence, taking the minimum value of
$2^{nH(d/n)}$ in the relevant range yields
$$
\sum_{i=0}^d {\binom{n}{i}} = 1 + .97 \cdot 2^{nH(d/n)} < .98 \cdot 2^{nH(d/n)}.
$$

Now, when $d>p$ (and $n \ge 41$) we have that
\beq
\sum_{i=1}^d {\binom{n}{i}}
& = & \sum_{i=1}^d \frac{1}{\sqrt{2\pi i(1-i/n)}} \cdot 2^{nH(i/n)} \\
& = & \sum_{i=1}^{\floor{p}} \frac{1}{\sqrt{2\pi i(1-i/n)}} \cdot 2^{nH(i/n)} +
      \sum_{i=\floor{p}+1}^d \frac{1}{\sqrt{2\pi i(1-i/n)}} \cdot 2^{nH(i/n)} \\
& < & \frac{1}{\sqrt{\pi}} (\frac{1}{\sqrt{\floor{p}}}+1) \cdot 2^{nH(d/n)} +
      \sum_{i=\floor{p}+1}^d \frac{1}{\sqrt{\pi i}} \cdot 2^{nH(i/n)} \\
& \le & \frac{1}{\sqrt{\pi}} (\frac{1}{\sqrt{\floor{p}}}+ 1
      +\frac{n/2-\floor{p}}{\sqrt{\floor{p}+1}} ) \cdot 2^{nH(d/n)} \\
& < & 0.9963 \cdot 2^{nH(d/n)},
\eeq
where the final inequality follows from taking the minimum value of $p$ (and noting that
$p$ achieves its minimum when $n$ achieves its minimum). Hence, taking the minimum value of
$2^{nH(d/n)}$ in the relevant range yields
$$
\sum_{i=0}^d {\binom{n}{i}} = 1 + .9963 \cdot 2^{nH(d/n)} < .9964 \cdot 2^{nH(d/n)}.
$$

It remains only to prove the above claim that
$\ln (\frac{n}{i} - 1) > \frac{1}{\sqrt{i}} $
for $1 \le i \le p$,
or equivalently that
$$g_n(i) = e^{-1/\sqrt{i}} (\frac{n}{i} - 1) > 1$$
in the desired range.
First note that $g_n(i)$ is strictly decreasing in this range, since
$g_n'(i) = \frac{-ni^{1/2}+n-i}{2i^{5/2}e^{1/\sqrt{i}}} < 0$, so we need only
demonstrate that $g_n(p)>1$. To this end, we consider the function
\beq
h_p(n)
= g_n(p)
= \left( \frac{2}{1-2/c\sqrt{n}} -1 \right) e^{-1/\sqrt{n/2 - \sqrt{n}/c}},
\eeq
and note that its derivative is negative in the relevant range:
\beq
h_p'(n)
& = & \frac{e^{-1/\sqrt{n/2-\sqrt{n}/c}}}{n^{3/2}(1-2/c\sqrt{n})^2}
      \left[ \frac{-2}{c} +
    \frac{1}{\sqrt{2}} \cdot
       \frac{1+2/c\sqrt{n}}{(1-2/c\sqrt{n})^{1/2}} \cdot
    (1-\frac{1}{c\sqrt{n}})
      \right] \\
& < & \frac{e^{-1/\sqrt{n/2-\sqrt{n}/c}}}{n^{3/2}(1-2/c\sqrt{n})^2}
      \left[ \frac{-2}{c} +
    \frac{1}{\sqrt{2}} \cdot
       \frac{1+2/c\sqrt{n}}{(1-2/c\sqrt{n})^{1/2}}
      \right] \\
& < & 0,
\eeq
where the final inequality follows by noting that the bracketed term is maximized
when $n$ takes its minimum value. It follows that $g_n(p) = h_p(n) \ge h_p(\infty) = 1$.

\enpf


\hide{
We will also need a
simple relation between the functions $H$ and $\beta$:
\belen
\label{lem:Hb}
For $d\geq1$ and $\gamma\in[0,1]$, we have
\beq
H\paren{\frac{ d}{\ceil{\beta(\gamma)d}}}
\leq
\log(2/(1+\gamma))
\hide{
H\paren{\frac{ d}{\floor{\beta(\gamma)d}}}
\leq
-\paren{1+
d\inv
}\log(\half+\tfrac{\gamma}{2})
}
,
\eeq
where $\beta(\cdot)$ is defined in (\ref{eq:betadef}).
\enlen
\bepf
Since $\beta(\gamma)\geq2$ implies $0<d/\ceil{\beta(\gamma)d}\leq\half$, we have
\beq
H\paren{\frac{ d}{\ceil{\beta(\gamma)d}}}
&\leq&
H\paren{\frac{ d}{{\beta(\gamma)d}}} \\
&=& H(1/\beta(\gamma)) \\
&=& \log(2/(1+\gamma))
\eeq
by definition of $\beta$.
\hide{
We note the elementary relation
\beq
x-\oo{k}\leq \frac{\floor{kx}}{k} \leq x,
\eeq
valid for all $x\in\R$ and $k>0$. It follows that
\beq
\frac{ d}{\floor{\beta(\gamma)d}}
\leq \oo{\beta-\oo{d}}
= \oo{\beta}+\oo{\beta(\beta d-1)}
\leq
\oo{\beta}+\oo{\beta d}
\leq 1
,
\eeq
since $\beta(\gamma)\geq2$.
Since $H$ is concave, we have
\beq
 H(x+\delta)\leq H(x)+\delta H(x')
\eeq
for all $x,\delta\geq0$ with $x+\delta\leq1$.
Additionally, we observe that
\beq
H'(x) &=& -\log(x)+\log(1-x) \\
&\leq& -\log(x)-\frac{1-x}{x}\log(1-x) \\
&=& \frac{H(x)}{x}
\eeq
for all $0\leq x\leq 1$.
Therefore,
\beq
H\paren{\frac{ d}{\floor{\beta(\gamma)d}}} &\leq &
H\paren{ \oo{\beta}+\oo{\beta d} } \\
&\leq& H\paren{\oo{\beta}} + \oo{\beta d}H'\paren{\oo{\beta}} \\
&\leq& H\paren{\oo{\beta}} + \oo{d}H\paren{\oo{\beta}} \\
&=&
-\paren{1+
d\inv
}\log(\half+\tfrac{\gamma}{2})
\eeq
by definition of $\beta$.
}
\enpf
}
}

Our extension to $k$-ary alphabets requires the corresponding analogue of
Lemma \ref{lem:binom}:

\belen
\label{lem:binom-k}
For $2 \le d \le \frac{k}{k+
1.6
} \cdot n$
and $n \ge 6$, we have
\beq
\sum_{i=0}^d \binom{n}{i}k^i &<& .94 \cdot 2^{nH(d/n) + d \log k}.
\eeq
\enlen

\bepf
First note that the derivative of $f(i)=2^{nH(i/n) + i \log k}$ is
$f'(i) = f(i) [\ln (\frac{n}{i} - 1) + \ln k]di$,
so $f(i)$ attains its maximum over the
range $0 \le i \le n$ at $i = \frac{k}{k+1} \cdot n$. Further note that
since
$i \le d \le \frac{k}{k+1.6} \cdot n
< \frac{k}{k+e^{1/\sqrt{\floor{n/2}+1}}} \cdot n$, we have that
$\ln (\frac{n}{i} - 1) + \ln k > \frac{1}{\sqrt{\floor{n/2}+1}}$.

We break up the analysis into two cases: When $d \le \frac{n}{2}$ we have
\beq
\sum_{i=0}^d \binom{n}{i}k^i
&<&
\binom{n}{d} 2k^d \\
&<& \frac{2 k^{d}}{\sqrt{\pi d}} 2^{nH(d/n)} \\
&=& \frac{2}{\sqrt{\pi d}} f(d) \\
&<& .8 \cdot f(d)
.
\eeq
When $d > \frac{n}{2}$ we have
\beq
\sum_{i=0}^d \binom{n}{i}k^i
&=&
\sum_{i=0}^{\floor{n/2}} \binom{n}{i}k^i
+ \sum_{i=\floor{n/2}+1}^d \binom{n}{n-i}k^i \\
&<&
\frac{2f(\floor{n/2})}{\sqrt{\pi \floor{n/2}}}
+ \sum_{i=\floor{n/2}+1}^d \frac{2^{nH(1-i/n)} k^i}{\sqrt{2\pi i(i/n)}} \\
&<&
\frac{2f(\floor{n/2})}{\sqrt{\pi \floor{n/2}}}
+ \frac{1}{\sqrt{\pi (\floor{n/2}+1)}} \sum_{i=\floor{n/2}+1}^d 2^{nH(i/n) + i \log k} \\
&<&
\frac{2f(\floor{n/2})}{\sqrt{\pi \floor{n/2}}}
+ \frac{f(d)}{\sqrt{\pi (\floor{n/2}+1)}}
+ \frac{1}{\sqrt{\pi (\floor{n/2}+1)}} \sum_{i=\floor{n/2}+1}^{d-1} f(i) \\
&<&
\frac{2f(\floor{n/2})}{\sqrt{\pi \floor{n/2}}}
+ \frac{f(d)}{\sqrt{\pi (\floor{n/2}+1)}}
+ \frac{1}{\sqrt{\pi (\floor{n/2}+1)}} \int_{\floor{n/2}+1}^{d} f(i) \\
&\le&
\frac{2f(\floor{n/2})}{\sqrt{\pi \floor{n/2}}}
+ \frac{f(d)}{\sqrt{\pi (\floor{n/2}+1)}}
+ \frac{1}{\sqrt{\pi}} \int_{\floor{n/2}+1}^{d} f(i) [\ln (\frac{n}{i} - 1) + \ln k]di \\
&=&
\frac{2f(\floor{n/2})}{\sqrt{\pi \floor{n/2}}}
+ \frac{f(d)}{\sqrt{\pi (\floor{n/2}+1)}}
+ \frac{f(d)}{\sqrt{\pi}}
- \frac{f(\floor{n/2}+1)}{\sqrt{\pi}} \\
&<& .94 f(d).
\eeq

\enpf

\section*{Acknowledgements}
We thank Noga Alon for helpful comments and references about Theorem \ref{thm:vclb}
and the anonymous referee for weeding out a number of inaccuracies.


\begin{thebibliography}{10}

\bibitem{alon97scalesensitive}
Noga Alon, Shai Ben-David, Nicol{\`o} Cesa-Bianchi, and David Haussler.
\newblock Scale-sensitive dimensions, uniform convergence, and learnability.
\newblock {\em Journal of the ACM}, 44(4):615--631, 1997.

\bibitem{dana-david-me-lev/SQ}
Dana Angluin, David Eisenstat, Leonid Kontorovich, and Lev Reyzin.
\newblock Lower bounds on learning random structures with statistical queries.
\newblock In \emph{ALT}, pages 194--208, 2010.

\bibitem{DBLP:conf/stoc/BlumFJKMR94}
Avrim Blum, Merrick~L. Furst, Jeffrey~C. Jackson, Michael~J. Kearns, Yishay
  Mansour, and Steven Rudich.
\newblock Weakly learning dnf and characterizing statistical query learning
  using fourier analysis.
\newblock In {\em STOC}, pages 253--262, 1994.

\bibitem{MR1072253}
Anselm Blumer, Andrzej Ehrenfeucht, David Haussler, and Manfred~K. Warmuth.
\newblock Learnability and the {V}apnik-{C}hervonenkis dimension.
\newblock {\em J. Assoc. Comput. Mach.}, 36(4):929--965, 1989.

\bibitem{Bshouty2009323}
Nader~H. Bshouty, Yi~Li, and Philip~M. Long.
\newblock Using the doubling dimension to analyze the generalization of
  learning algorithms.
\newblock {\em Journal of Computer and System Sciences}, 75(6):323 -- 335,
  2009.

\bibitem{MR512411}
R.~M. Dudley.
\newblock Central limit theorems for empirical measures.
\newblock {\em Ann. Probab.}, 6(6):899--929 (1979), 1978.

\bibitem{feldman09}
Vitaly Feldman.
\newblock A complete characterization of statistical query learning with
  applications to evolvability.
\newblock In {\em Symposium on Foundations of Computer Science (FOCS)}, 2009.

\bibitem{1121738}
J.~Flum and M.~Grohe.
\newblock {\em Parameterized Complexity Theory (Texts in Theoretical Computer
  Science. An EATCS Series)}.
\newblock Springer-Verlag New York, Inc., Secaucus, NJ, USA, 2006.

\bibitem{Gilbert52}
E.~N. Gilbert.
\newblock A comparison of signalling alphabets.
\newblock {\em Bell Syst Tech J}, 31:504--522, 1952.

\bibitem{MR757767}
Evarist Gin{\'e} and Joel Zinn.
\newblock Some limit theorems for empirical processes.
\newblock {\em Ann. Probab.}, 12(4):929--998, 1984.
\newblock With discussion.

\bibitem{DBLP:conf/focs/Haussler89}
David Haussler.
\newblock Generalizing the PAC model: Sample size bounds from metric
  dimension-based uniform convergence results.
\newblock In \emph{30th Annual Symposium on Foundations of Computer Science},
  pages 40--45, 1989.

\bibitem{MR1313896}
David Haussler.
\newblock Sphere packing numbers for subsets of the {B}oolean {$n$}-cube with
  bounded {V}apnik-{C}hervonenkis dimension.
\newblock {\em J. Combin. Theory Ser. A}, 69(2):217--232, 1995.

\bibitem{209962}
David Haussler and Philip~M. Long.
\newblock A generalization of sauer's lemma.
\newblock {\em J. Comb. Theory Ser. A}, 71(2):219--240, 1995.

\bibitem{MR1965359}
S.~Mendelson and R.~Vershynin.
\newblock Entropy and the combinatorial dimension.
\newblock {\em Invent. Math.}, 152(1):37--55, 2003.

\bibitem{DBLP:journals/ml/Natarajan89}
B.~K. Natarajan.
\newblock On learning sets and functions.
\newblock \emph{Machine Learning}, 4: 67--97, 1989.

\bibitem{pollard84}
David Pollard.
\newblock \emph{Convergence of Stochastic Processes}.
\newblock Springer-Verlag, 1984.

\bibitem{MR1089429}
David Pollard.
\newblock {\em Empirical processes: theory and applications}.
\newblock NSF-CBMS Regional Conference Series in Probability and Statistics, 2.
  Institute of Mathematical Statistics, Hayward, CA, 1990.

\bibitem{RothSeroussi}
Ron~M. Roth and Gadiel Seroussi.
\newblock Bounds for binary codes with narrow distance distributions.
\newblock {\em IEEE Transactions on Information Theory}, 53(8):2760--2768,
  2007.

\bibitem{MR0307902}
Norbert Sauer.
\newblock On the density of families of sets.
\newblock {\em J. Combinatorial Theory Ser. A}, 13:145--147, 1972.

\bibitem{MR905334}
Michel Talagrand.
\newblock Donsker classes and random geometry.
\newblock {\em Ann. Probab.}, 15(4):1327--1338, 1987.

\bibitem{talagrand87}
Michel Talagrand.
\newblock {The Glivenko-Cantelli Problem}.
\newblock {\em Ann. Probab.}, 15(3):837--870, 1987.

\bibitem{MR945108}
Michel Talagrand.
\newblock Donsker classes of sets.
\newblock {\em Probab. Theory Related Fields}, 78(2):169--191, 1988.

\bibitem{vapnik82}
Vladimir~N. Vapnik.
\newblock {\em Estimation of dependences based on empirical data}.
\newblock Springer-Verlag, 1982.

\bibitem{Varshamov57}
R.~R. Varshamov.
\newblock Estimate of the number of signals in error correcting codes.
\newblock {\em Dokl Acad Nauk SSSR}, 117:739--741, 1957.

\end{thebibliography}
\end{document}